\documentclass[reqno,10pt,centertags,draft]{amsart}
\usepackage{amsmath,amsthm,amscd,amssymb,latexsym,upref}
\date{\today}
\newcommand{\bbN}{{\mathbb{N}}}

\newcommand{\bbC}{{\mathbb{C}}}

\renewcommand{\Re}{\text{\rm Re}}

\def\sgn{{\mathrm{sgn}\,}}
\renewcommand{\Im}{\text{\rm Im}}

\allowdisplaybreaks
\numberwithin{equation}{section}
\newtheorem{theorem}{Theorem}[section]

\newtheorem{lemma}[theorem]{Lemma}

\theoremstyle{definition}

\newtheorem{problem}[theorem]{Problem}

\begin{document}

\title[]{Uniform approximation of $\sgn(x)$ by rational functions with
prescribed poles}
\author[]{F. Peherstorfer,  and P. Yuditskii}

\thanks{
Partially supported by  
the Austrian
Founds FWF, project number: P16390--N04 and
{\it Marie Curie International
Fellowship} within the 6-th European Community Framework Programme,
 Contract MIF1-CT-2005-006966}
\dedicatory{Dedicated to the memory of B. Ya. Levin.}

\date{\today}

\maketitle

\begin{abstract} For $a\in (0,1)$  let $L^k_m(a)$ be the error of the
best approximation of the function $\sgn(x)$  
on the two symmetric intervals 
$[-1,-a]\cup[a,1]$
by  rational functions
with  the only possible poles of degree $2k-1$ at the origin and of
$2m-1$ at  infinity.  Then the following limit exists
\begin{equation}\label{0.1}
\lim_{m\to \infty}L^k_m(a)\left(\frac{1+a}{1-a}\right)^{m-\frac{1}{2}}
(2m-1)^{k+\frac{1}{2}}=\frac
2 \pi\left(\frac{1-a^2}{2a}\right)^{k+\frac{1}{2}}
\Gamma\left(k+\frac 1 2\right).
\end{equation}
%For the constants $\mu(p)>0$ we have representations in terms of the
%special explicitly given conformal mappings. Note that for 
%unweighted approximation error  $E^*_n(p,a)$ and for all
%$a$ we have
% $E_n(p,a) C_1(a)\le E^*_n(p,a)\le C_2(a) E_n(p,a)$, and moreover
% $\lim_{a\to 1}\lim_{n\to\infty}\frac{E^*_n(p,a)}{E_n(p,a)}=1$.

\end{abstract}

\section{Introduction}
This is the second step (for the first one see \cite{EYU})
on the way to understand better the difficulties that up  to now are not
allow to find the Bernstein constant. Recall that Sergey Natanovich
Bernstein found  \cite{Bern1,Bern2}
 that for the error $E_n(p)$ of the best uniform approximation
 of $|x|^p$, $p$  not an even integer, on $[-1,1]$
by polynomials of degree $n$ the following limit exists:
$$\lim_{m\to\infty}n^p E_n(p)=\mu(p)>0.$$
For $p=1$ this result  was obtained by Bernstein in 1914,
and he posed the question, whether one can express
$\mu(1)$ in terms of some known transcendental functions.
This question is still open. 

Actually, we solve here a  problem on asymptotics
of the best approximation of $\sgn(x)$ 
on the union of two intervals $[-1,-a]\cup[a,1]$
by rational functions.
In 1877, E. I. Zolotarev \cite{Zol,Akh}
found an explicit expression, in terms of elliptic functions,
of the rational function
of given degree which is uniformly closest to $\sgn(x)$
on this set.
This result was subject of many generalizations, and it has
applications in electric engineering. 
In the Zolotarev's case the position of poles of the rational
function is free,
a natural question is to find the best approximation
when the poles and their multiplicities are   fixed.
In \cite{EYU} A. Eremenko and the second author solved the polynomial
case. Here we allow the rational function to have one more pole
in $(-a,a)$, more precisely, we  allow two
poles --- one at infinity and one in the origin.

Thus the problem is: 

\begin{problem}\label{p1.1} For $k,m\in\bbN$, 
find the best
approximation of the function $\sgn(x)$,
 $|x|\in [a,1]$, by functions of the form
$$
f(x)=\frac{a_{-(2k-1)}}{x^{2k-1}}+...+{a_{2m-1}}{x^{2m-1}}
$$
and   the
approximation error $L^k_m(a)$.
\end{problem}

  One can be interested in many different
asymptotics for $L^k_m(a)$ when $m$ or $k$ or both 
 go to infinity in a certain
prescribed way. In this work we concentrate
on the case when $k$ is fixed and $m\to \infty$.
Note, however, that due to the evident symmetry
$L^k_m(a)=L^m_k(a)$ and a bit less evident \eqref{6.2.15}
we have simultaneously asymptotic for $k\to \infty$, $m$ is fixed
and $k\to\infty$, $m\to \infty$ so that $k=m$.

As it appears the tricks which are used in \cite{EYU} to find precise
asymptotic
 work in this general case (so we have a method in hands):
\begin{itemize}
\item[1.] For each particular $k$ and $m$ we reveal the structure of the
extremal function by representing it with the help of an explicitly
given conformal mapping.
\item[2.] The system of conformal mappings ($k$ is fixed, $m$ is a
parameter) converges (in the Caratheodory sense) after an appropriate
renormalization. The limit map does not depend on $a$, thus we obtain
asymptotics for $L^k_m(a)$ in terms of $a$--depending parameters,
that we use for renormalization, (an explicit formula)
and a $k$--depending constant say $Y_k$, which is a
certain characteristic of this final conformal map (kind of capacity). 
\end{itemize}

\noindent
Of course, it is very tempting to guess $Y_k$ directly from the 
given explicitly conformal map. It might be that we have here
special functions that are given in such form that we are unable
to recognize them. In any case, we would consider this way of
finding $Y_k$ as a very interesting open problem. However we are
able to find $Y_k$ using the third step below
of our strategy. Problem \ref{p1.1}
in an evident way is equivalent to
\begin{problem}\label{p1.2}
For $p=2k-1$ and $n=2(k+m-1)$, find the best weighted
{\it polynomial} approximation and the minimal deviation
\begin{equation}
E^*_n(p,a)=\inf_{\{P:\deg P\le n\}}\sup_{|x|\in[a,1]}
\left\vert\frac{|x|^{p}-P(x)}{x^{p}}\right\vert.
\end{equation}
\end{problem}
\noindent
Thus we have
$E^*_n(p,a)=L_m^k(a)$. Note that Bernstein himself solved the unweighted
\begin{problem}\label{p1.3}
For a fixed non even $p$, find asymptotics for
 the minimal deviation
\begin{equation}
E_n(p,a)=\inf_{\{P:\deg P\le n\}}\sup_{|x|\in[a,1]}
\left\vert{|x|^{p}-P(x)}\right\vert,
\end{equation}
when $n$ goes to infinity through the even integers.
\end{problem}
\begin{itemize}
\item[3.]
Due to the evident relation
$$
\lim_{a\to 1}\lim_{n\to \infty}\frac{E^*_n(p,a)}{E_n(p,a)}=1,
$$
we can recalculate the constant in Problem \ref{p1.3}
to the constant related to  Problem \ref{p1.2} and thus to get
explicitly $e^{Y_k}=\frac{\Gamma\left(k+\frac 1 2
\right)}{2^{k+\frac 1 2}\pi}$.
\end{itemize}

This interplay between Problems \ref{p1.2} and \ref{p1.3} indicates
that most likely one can find our asymptotic formula \eqref{0.1}
by using original Bernstein's method, though up to the last step
our consideration are very direct and simple. However we can go in the
opposite direction. In particular in this work we show that the extremal
polynomials of Problem \ref{p1.3}, at least for $p=1$, also have 
special representations in terms of conformal mappings. The boundary
of the corresponding domains are not so explicit as in Problem
\ref{p1.1}, they
are described in terms of certain functional equations involving unknown
function, its Hilbert transform and independent variable
\eqref{6.2.18}. Precise constants that characterize these equations
(counterparts of the constants $Y_k$),
related to the conformal mappings and their asymptotics leave enough
space for the hope that for $a=0$ one also would be able to characterize
very similar equations in terms of classical constants.

\bigskip

\noindent{\bf Acknowledgment.} 
We are thankful to Alex Eremenko for friendly conversations during the
writing of this paper.

\section{Special Functions}

In this section we introduce certain special conformal mappings that we
need in what follows. They are marked by a natural parameter $k$, but in
this section $k$ can be just real, $k>1/2$.

For given  $k$, consider the domain 
\begin{equation}
\Pi_k=\bbC_+\setminus\{w: \Re w=-\log t,\ |\Im w-k\pi|\le
\arccos t, \ t\in (0,1]\}
\end{equation}
Define the conformal map
$$
H_k:\bbC_+\to \Pi_k
$$
normalized by $H_k(0)=\infty_1$,
$H_k(\infty)=\infty_2$ (on the boundary we have two infinite points
that we denote respectively $\infty_1,\infty_2$), and moreover
$$
H_k(\zeta)=\zeta+..., \quad \zeta\to \infty,
$$
(that is the leading coefficient is fixed).
By $D_k$ we denote the positive number such that $H_k(-D_k)=0$.

Note that for $H_k$ we have the following integral representation
\begin{equation}
H_k(\zeta)=\zeta+D_k+\int_0^\infty
\left(\frac{1}{t-\zeta}-\frac{1}{t+D_k}\right)\rho_k(t) dt,
\end{equation}
where $\rho_k(t)=\frac{1}{\pi}\Im H_k(t)$. Evidently
$\rho_k(t)\to k+\frac 1 2$, $t\to +\infty$.
\begin{lemma} The function $H_k$ possesses the asymptotic
\begin{equation}\label{2.3}
\lim_{\zeta\to-\infty}
\left\{
H_k(\zeta)-\zeta+\left(k+\frac 1 2\right)\log(-\zeta)
\right\}=Y_k,
\end{equation}
where
\begin{equation}
Y_k:=D_k+
\left(k+\frac 1 2\right)\log D_k
-\int_0^\infty
\frac{\rho_k(t)-\left(k+\frac 1 2\right)}{t+D_k}
dt.
\end{equation}

\end{lemma}
\begin{proof}
Since
\begin{equation}
\int_0^\infty
\left(\frac{1}{t-\zeta}-\frac{1}{t+D_k}\right)
\left(\rho_k(t)-\left(k+\frac 1 2\right)\right)
dt\to
-\int_0^\infty
\frac{\rho_k(t)-\left(k+\frac 1 2\right)}{t+D_k}
dt,
\end{equation}
and
\begin{equation}
\left(k+\frac 1 2\right)\int_0^\infty
\left(\frac{1}{t-\zeta}-\frac{1}{t+D_k}\right)
dt=-\left(k+\frac 1 2\right)(\log(-\zeta)-\log D_k)
\end{equation}
we get \eqref{2.3}.
\end{proof}

Finally note that $Y_k$, as it was defined here, has sense for all real
$k>\frac 1 2$. As it is shown in Sect. 5 for an integer $k$ 
we have
$$
Y_k=\log{\Gamma\left(k+\frac 1 2
\right)} -\left(k+\frac 1 2\right)\log 2-\log\pi.
$$
We do not know if these values coincide for non integers $k$.

\section{Extremal problem}
Problems \ref{p1.1} and \ref{p1.2} are related in a trivial way.
Recall,
for $p=2k-1$ and $n=2(k+m-1)$, we have
\begin{equation}
E^*_n(p,a)=L_m^k(a)=\inf_{\{P:\deg P\le 2(m+k-1)\}}\sup_{|x|\in[a,1]}
\left\vert\frac{|x|^{2k-1}-P(x)}{x^{2k-1}}\right\vert,
\end{equation}
where $a\in (0,1)$, $k,m\in \bbN$.
Evidently, $L_m^k(a)$ can be rewritten in the terms of the best
approximation of the function $\sgn(x)$ by functions of the form
$$
f(x)=\frac{a_{-(2k-1)}}{x^{2k-1}}+...+{a_{2m-1}}{x^{2m-1}}.
$$
Also, it is trivial that the extremal polynomial is even in the first
case and the extremal function $f=f(x;k,m;a)$ is odd.

For a parameter $B>0$ and $k,m\in \bbN$, $\Omega_m^k(B)$ denotes the
subdomain of the half strip
$$
\{w=u+iv: v>0, \ 0<u<(k+m)\pi\}
$$
that we obtain by deleting the subregion
\begin{equation}\label{8aug}
\{w=u+iv: |u-\pi k|\le\arccos\left(\frac{\cosh B}{\cosh v}\right),
\ v\ge B\}.
\end{equation}
Let $\phi(z)=\phi(z;k,m;B)$ be the conformal map of the first quadrant
onto $\Omega_m^k(B)$ such that
$\phi(0)=\infty_1$, $\phi(1)=(k+m)\pi$, $\phi(\infty)=\infty_2$.
Let $a=\phi^{-1}(0)$. Then $a$ is a continuous strictly increasing
function of $B$, moreover $\lim_{B\to 0} a(B)=0$ and
$\lim_{B\to\infty}a(B)=1$. Thus we may consider the inverse function
$B(a)=B^k_m(a)$, $a\in(0,1)$.
\begin{theorem} The error of the best approximation is
\begin{equation}\label{3.2}
L_m^k(a)=\frac{1}{\cosh B^k_m(a)}
\end{equation}
and the extremal function is of the form
$$
f(x;k,m;a)=1-(-1)^k L_m^k(a)\cos\phi(x;k,m;B(a)), \quad x>0.
$$
\end{theorem}

\begin{proof}
Basically the proof is the same as in \cite{EYU}. A comparably important
difference is as follows. We have to note and prove that on the imaginary
axis the extremal function has precisely one zero (there are no critical
points and the behavior at
$i0$ and at
$i\infty$ is evident). At this point $\phi=k\pi+iB$ and we have
\eqref{3.2}.
\end{proof}

\section{Asymptotics}
\begin{theorem} The following limit exists
\begin{equation}\label{4.1}
\begin{split}
\lim_{m\to\infty}\left\{B^k_m(a)-\left(m-\frac 1 2\right)
\log\frac{1+a}{1-a}-\left(k+\frac 1 2\right)\log (2m-1)
\right\}\\=\left(k+\frac 1 2\right)\log{\frac{a}{1-a^2}}-Y_k.
\end{split}
\end{equation}
\end{theorem}

\begin{proof}
As in \cite{EYU} we use the symmetry principle and make a convenient
changes of variable to have a conformal map $\Phi_m(Z)=\Phi(Z;k,m;B)$
of
the upper plane in  the region 
$$
i(\Omega_m^k(B)\cup\overline{\Omega_m^k(B)})\cup(0,i\pi(m+k)).
$$
This conformal map has the following boundary correspondence
$$
\Phi_m:(-C_m,-A_m,0,A_m,C_m)\to
(-\infty_2,-\infty_1,0,\infty_1,\infty_2),
$$
here $A_m=aC_m$ and $C_m$ will be chosen a bit later.

For $\Phi_m$ we have the following integral representation
$$
\Phi_m(Z)=\left(m-\frac 1 2\right)\log\frac{1+\frac{Z}{C_m}}
{1-\frac{Z}{C_m}}+\int_{A_m}^\infty
\left[\frac{1}{X-Z}-\frac{1}{X+Z}\right]v_m(X)\,dX,
$$
where
\begin{equation}
v_m(X)=
\begin{cases}
\frac 1{\pi}\Im \Phi_m(X),& A_m\le X\le C_m \\
k+\frac 1 2,& X> C_m
\end{cases}
\end{equation}
Put now
$$
H_m^k(\zeta)=\Phi_m(Z)-B_m,\quad Z=A_m+\zeta,
$$
then
$$
H^k_m(\zeta)=\left(m-\frac 1 2\right)\log\frac{1+a+\frac{\zeta}{C_m}}
{1-a-\frac{\zeta}{C_m}}+\int_{0}^\infty
\left[\frac{1}{t-\zeta}-\frac{1}{t+A_m+\zeta}\right]\hat v_m(t)\,dt
-B_m,
$$
where $\hat v_m(t)=v_m(t+A_m)$. Let us rewrite $H^k_m$ in the form that
is close to the integral representation of $H_k$:
\begin{equation}\label{4.2}
\begin{split}
H^k_m(\zeta)=&\left(m-\frac 1 2\right)\log\frac{1+\frac{\zeta}{C_m(1+a)}}
{1-\frac{\zeta}{C_m(1-a)}}+D_k+
\int_{0}^\infty
\left[\frac{1}{t-\zeta}-\frac{1}{t+D_k}\right]\hat v_m(t)\,dt\\
+&
\left(m-\frac 1
2\right)\log\frac{1+a}{1-a}-D_k
+\int_{0}^\infty
\left[\frac{1}{t+D_k}-\frac{1}{t+A_m+\zeta}\right]\hat v_m(t)\,dt
-B_m
\end{split}
\end{equation}

Now,  we put
$$
C_m=\frac{2m-1}{1-a^2}.
$$
In this case the first line in \eqref{4.2} converges to $H_k(\zeta)$.
Since
\begin{equation}\label{}
\begin{split}
\lim_{m\to\infty}\int_{0}^\infty
\left[\frac{1}{t+D_k}-\frac{1}{t+A_m+\zeta}\right]
\left(\hat v_m(t)-\left(k+\frac 1 2\right)
\right)\,dt\\
=\int_{0}^\infty
\frac{\rho_k(t)-\left(k+\frac 1 2\right)}{t+D_k}dt
\end{split}
\end{equation}
and
\begin{equation}\label{}
\int_{0}^\infty
\left[\frac{1}{t+D_k}-\frac{1}{t+A_m+\zeta}\right]
\,dt
=\log\frac{A_m}{D_k}+\log\left(1+\frac{\zeta}{A_m}\right)
\end{equation}
we have from the second line in \eqref{4.2} that
\begin{equation}\label{}
\begin{split}
\lim_{m\to\infty}\left\{B_m-\left(m-\frac 1
2\right)\log\frac{1+a}{1-a}-\left(k+\frac 1 2\right)\log{A_m}\right\}\\
=-D_k-\left(k+\frac 1 2\right)\log D_k+\int_{0}^\infty
\frac{\rho_k(t)-\left(k+\frac 1
2\right)}{t+D_k}dt=-Y_k.
\end{split}
\end{equation}
Thus we get \eqref{4.1}. In order to prove \eqref{0.1}
we have to find
the constant $2e^{Y_k}$.

\end{proof}

\section{The constant}

From the point of view of the best weighted polynomial approximation of
the function $|x|^p$ (see Sect. 3)  our current result has the form
\begin{equation}
\lim_{m\to \infty}\left(\frac{1+a}{1-a}\right)^{\frac{n}{2}+1}
n^{\frac{p}{2}+1}E^*_n(p,a)=\left(\frac{(1+a)^2}{a}\right)^{\frac{p}{2}+1}
c(p).
\end{equation}
On the other hand for the uniform approximation
of $|x|^p$ (details see in Appendix~1)
\begin{equation}
\lim_{m\to \infty}\left(\frac{1+a}{1-a}\right)^{\frac{n}{2}+1}
n^{\frac{p}{2}+1}E_n(p,a)={2}^{\frac{p}{2}+1}{a}^{\frac{p}{2}-1}
\frac{(1+a)^2}{2\left|\Gamma\left(-\frac p 2\right)\right|}.
\end{equation}
Since
$$
\lim_{a\to 1}\lim_{n\to \infty}\frac{E^*_n(p,a)}{E_n(p,a)}=1,
$$
we obtain
$$
{c(p)2^{\frac p 2 +1}}\left|\Gamma\left(-\frac p 2\right)\right|=2.
$$
Using $\left|\Gamma\left(-\frac p 2\right)\right|
\Gamma\left(\frac p 2+1\right)=\pi$, we have
$$
c(p)=\frac 2 \pi 2^{-\frac p 2 -1}\Gamma\left(\frac p
2+1\right).
$$
This finishes the proof of \eqref{0.1}.

\section{Case $m=k$, $m\to \infty$}
It is quite evident that the final configuration of the
conformal mapping in this case should be just a symmetrization of
the map that we had in the case $k=0$,
$m\to \infty$. However it's even much simpler to make this reduction by a
suitable change of variable.\
First we put $a=\alpha^2$, then  $x\in[a,1]$ means
$y=\frac{x}{\alpha}\in[\alpha,\alpha^{-1}]$ and we have one more symmetry
$y\mapsto 1/y$. Therefore the extremal function is symmetric and
possesses the representation
\begin{equation}
\tilde f(y;m,m):=f(x;m,m;a)=P_{2m-1}\left(\frac{y+y^{-1}}
{\alpha+\alpha^{-1}}\right),
\end{equation}
where $P_{2m-1}(t)$ is the best polynomial approximation of
$\sgn(t)$ on 
$\left[-1,-\frac{2\alpha}{1+\alpha^{2}}\right]
\cup\left[\frac{2\alpha}{1+\alpha^{2}},1\right]$.
Thus we just have
\begin{equation}\label{6.2.15}
L^m_m(a)=L_m^0\left(\frac{2\sqrt{a}}{1+a}\right),
\end{equation}
and
\begin{equation}\label{6.2}
\lim_{m\to
\infty}L^m_m(a)\left(\frac{1+\sqrt{a}}
{1-\sqrt{a}}\right)^{2m-{1}}
(2m-1)^{\frac{1}{2}}=
\frac{1-a}{\sqrt{\pi\sqrt{a}(1+a)}}.
\end{equation}

\section{Unweighted extremal polynomial via conformal mapping}
Let $P_m(z,a)$ be the best uniform (unweighted) approximation
of $|x|$ by polynomials of degree not more than $2m$ on
two intervals $[-1,-a]\cup[a,1]$
and $L=L_m(a)$ be the approximation error.

In this
section we prove
\begin{theorem}
There is a curve $\gamma=\gamma_m(a)$ inside the half--strip
\begin{equation}\label{6.1.18}
\{w=u+iv: u\in (0,(m+1)\pi),\ v>0\}
\end{equation}
such that the extremal polynomial possesses the representation
$$
P_m(z,a)=z+L\cos\phi_m(z,a)
$$
where $\phi_m(z,a)$ is the conformal map of the first quadrant onto the
region in the half strip \eqref{6.1.18} bounded on the left
by $\gamma_m(a)$, which is normalized by $\phi_m(a,a)=0$,
$\phi_m(1,a)=(m+1)\pi$ and $\phi_m(\infty,a)=\infty$. 
Moreover, the curve $\gamma$ is the image of the imaginary half--axis
under this conformal map that satisfies the following functional equation
\begin{equation}\label{6.2.18}
\gamma_m(a)=\{u+iv=\phi_m(iy,a):\ L\sin u(y)\sinh v(y)=y,\ y>0\}.
\end{equation}
\end{theorem}

\begin{proof}
First we clarify the shape of the extremal polynomial. In particular,
we prove that $P_m(0,a)>L$. On the way we show the fact that is probably
interesting on its own:
 $P'_m(x,a)$ looks pretty similar to  the polynomial of the best
approximation of $\sgn(x)$ on two symmetric intervals \cite{EYU}, with
the only difference that the deviations in area should be equal, instead
of the maximum modulus. 
However it can be shown
that  $P'_m(x,a)$ is not  the best
$L^1$ approximation of $\sgn(x)$.

Due to the symmetry of  $P_m(x,a)$, we can use the Chebyshev theorem with
respect to the best approximation of $\sqrt{x}$ on $[a^2,1]$ by
polynomials of degree
$m$. It gives us that
$P_m(z,a)$ has $m+2$ points $\{x_j\}$  on interval
$[a,1]$ where $P_m(x_j,a)=x_j\pm L$ (the right half of the Chebyshev
set in this case). Moreover,
$x_0=a$ and $x_{m+1}=1$.  At all other points, in addition, we have
$P'_m(x_j,a)=1$, $1\le j\le m$. Between each two of them 
we have a point $y_j$, where $P''_m(y_j)=0$. Therefore we 
obtain $2(m-1)$
zeros of the second derivative in $(-1,-a)\cup(a,1)$ and this is
precisely its degree. Thus there is {\it no other critical points} of
$P'_m(z,a)$, in particular, in $(-a,a)$ and on imaginary axis.

From the first consequence, we conclude that on $(-a,a)$ the $P'_m(z,a)$
increases. That is on $(a,x_1)$ the graph of $P_m(z,a)$ is under the line
$x\pm L$, depending on the value $P_m(a,a)$, that, recall,
should be $a+L$ or $a-L$. Therefore, it
is under the line $x+L$ and
$P_m(a,a)-a=L$, $P_m(x_1,a)-x_1=-L$.  
Continuing in this way we get values of $P_m(x_j,a)$ at all other points
$x_j$ by alternance principle.
Note that as byproduct we get
$$
\int_{x_{i-1}}^{x_i} |P'_m(x,a)-1|\,dx=2L
$$
for all $1\le i\le m+1$.

From the second consequence we have that $\Im P'_m(iy)\ge 0$ on the
imaginary axis, that is $P_m(iy,a)$, being real, decreases with $y$,
starting from $P_m(0,a)>L$ to $-\infty$. From this remark and the
argument principle we deduce that the equation
\begin{equation}\label{6.3.18}
P_m(z,a)-z=tL
\end{equation}
has no solution in the open first quarter for all $t\in (-1,1)$.

Indeed, since $P_m(z,a)-z$ alternate between $\pm L$
in the interval [a,1],  \eqref{6.3.18} has $m+1$ solutions, which we
denote by $x_j(t)$. Consider now the contour that runs on the positive
real axis till $x_j(t)-\epsilon$, then it goes around $x_j(t)$ on the
half--circle of the radius $\epsilon$
clockwise. After the last of $x_j$'s we continue to go 
along the contour
till the big
positive
$R$. Next piece of the contour is a quarter--circle till imaginary axis.
Finally, from
$iR$ we go back to the origin. On each half--circle of the radius
$\epsilon$ the argument of the function changes by $-\pi$.  On the
quarter circle it changes by about $\deg P_m(z,a)\times\frac{\pi}2=m\pi$.
On the imaginary axis we have $\Re(P_m(iy,a)-iy)= P_m(iy,a)$ and
$\Im(P_m(iy,a)-iy)= -y$. Since $P_m(iy,a)$ decreases and much faster than
$-y$ (degree of $P_m$ is at least two), the change of the argument
on the last piece of the contour is about $\pi$. Thus the whole change
is $-(m+1)\pi+m\pi+\pi=0$. Since the function has no poles, it has no
zeros in the region. 

Thus $\arccos\frac{P_m(z,a)-z} L$ is well define in the quarter--plane.
We finish the proof by  inspection of the boundary
correspondence. 
\end{proof}
Note two facts: the curve \eqref{6.2.18} has the asymptote
$u\to \pi$, $v\to+\infty$ ($y\to+\infty$) and we have uniqueness
of the solution of the functional equation  \eqref{6.2.18}
due to uniqueness of the extremal polynomial.

\section{Appendix 1}
From \cite{Akha}, problem 42:
\begin{equation}\label{6.1}
E_l\left[\frac{1}{(b+x)^s}\right]
\sim
\frac{l^{s-1}}{|\Gamma(s)|}
\frac{(b-\sqrt{b^2-1})^l} 
{(b^2-1)^{\frac{s+1}{2}}}\quad (b>1, \ s\not= 0),
\end{equation}
where $E_l[f(x)]$ is the error of the approximation of 
the function $f(x)$ on the interval $[-1,1]$ by polynomials of degree
not more than $l$.

We change the variable 
$$
y=\frac{b+x}{b+1}
$$
and put $a^2=\frac{b-1}{b+1}$. Then we have
$$
\inf_{P:\deg P\le
l}\max_{y\in [a^2,1]}|y^{-s}-P(y)|
=(1+b)^s E_l\left[\frac{1}{(b+x)^s}\right].
$$
That is
\begin{equation}\label{6.2}
E_{2l}(-2s,a)=(1+b)^s E_l\left[\frac{1}{(b+x)^s}\right].
\end{equation}

Note that
$$
b=\frac{1+a^2}{1-a^2},\quad 
b^2-1=\frac{4a^2}{(1-a^2)^2},
$$
and therefore
$$
\sqrt{b^2-1}=\frac{2a}{1-a^2},
\quad 
b-\sqrt{b^2-1}=\frac{1-a}{1+a}.
$$
Thus from \eqref{6.1} and \eqref{6.2}
we get
\begin{equation*}
\begin{split}
E_{2l}(-2s,a)\sim &
\left(\frac{2} 
{1-a^2}\right)^{s}
\frac{l^{s-1}}{|\Gamma(s)|}
\left(\frac{1-a} 
{1+a}\right)^l
\left(\frac{1-a^2} 
{2a}\right)^{s+1}\\
= &
a^{-s}
\frac{l^{s-1}}{|\Gamma(s)|}
\left(\frac{1-a} 
{1+a}\right)^l
\left(
\frac{1-a^2} {2a}
\right)\\
=&
a^{-s-1}
\frac{l^{s-1}}{|\Gamma(s)|}
\left(\frac{1-a} 
{1+a}\right)^{l+1}
\frac{(1+a)^2} 2.
\end{split}
\end{equation*}

\section{Appendix 2}

Here we present a "solvable model" for the problem under
consideration: we replace the comparably complicated
configuration \eqref {8aug}, that we remove from the strip,
by just two slits. We used this model on the first step of rough
understanding  a form of the asymptotic and  it might be useful for a
reader, in particular,
it contains the hint that in non--model case the asymptotic of
$L^{qm}_m(a)$ for $m\to \infty$  can also be found for an arbitrary 
$q\in
\bbN$ fixed, see \eqref{last}.

For $B>0$, consider the conformal map $w=\phi(z)$ of the upper 
half plane $\bbC_+$
onto  the strip
\begin{equation}
\Pi=\{w: 0<\Im w<(k+m)\pi\}
\end{equation}
with the cut
\begin{equation}
\gamma_B=\{w: \Im w=k\pi, \vert\Re w\vert\ge B\},
\end{equation}
under the normalizations
\begin{equation}
\phi(0)=0,\quad \phi(\pm 1)=\pm\infty_2,
\end{equation}
where $\infty_2$ denote the point on the boundary of  the domain when
we go to infinity on the level $k\pi<\Im w<(k+m)\pi$. By $\infty_1$ we
denote the point on the boundary  that corresponds to the level $0<\Im
w<k\pi$. Put
$a=\phi^{(-1)}(+\infty_1)$ (therefore $-a=\phi^{(-1)}(-\infty_1)$).

Let us find the precise formula for this map as well as the relation
between $a$ and $B$. We have
\begin{equation}
\begin{split}
\phi(z)=&k\int_a^\infty\left(\frac 1{x-z}-\frac 1{x+z}\right)dx+
m\int_1^\infty\left(\frac 1{x-z}-\frac 1{x+z}\right)dx\\
+&k\left.\log\frac{x-z}{x+z}\right\vert_a^\infty+
m\left.\log\frac{x-z}{x+z}\right\vert_1^\infty\\
=&k\log\frac{a+z}{a-z}+m\log\frac{1+z}{1-z}.
\end{split}
\end{equation}

Further, 
for $a<x<1$ we have
\begin{equation}
\Re \phi(x)=k\log\frac{x+a}{x-a}+m\log\frac{1+x}{1-x}
\end{equation}
and $B$ corresponds to the critical value of this function on the given
interval. For the critical point $c$ we have
\begin{equation}
\left(\Re \phi\right)'(c)=-\frac{2ka}{c^2-a^2}+\frac{2m}{1-c^2}=0.
\end{equation}
Therefore
\begin{equation}
c=\sqrt{\frac{ma^2+ka}{m+ka}}
\end{equation}
and
\begin{equation}
B=k\log\frac{c+a}{c-a}+m\log\frac{1+c}{1-c}.
\end{equation}

Let us mention that the relation between $a$ and $B$ is monotonic, and
$a$ runs from $0$ to $1$ as $B$ runs from $0$ to $\infty$.

As the next step we calculate the asymptotic behavior of $B$ for the
fixed $a$ as $m\to\infty$. First we write the asymptotic for $c$
\begin{equation}
c=\sqrt{\frac{ma^2+ka}{m+ka}}=a+\frac{k}{2m}(1-a^2)+...
\end{equation}
Therefore
\begin{equation}
\begin{split}
B=&k\log\left(2a+\frac{k}{2m}(1-a^2)+...\right)-k
\log\left(\frac{k}{2m}(1-a^2)+...\right)\\
+&m\log\frac{1+a+\frac{k}{2m}(1-a^2)+...}{1-a-\frac{k}{2m}(1-a^2)+...}\\
=&k\log\frac{2a}{1-a^2}+k
\log\frac{2m}{k}+...\\
+&m\log\frac{1+a}{1-a}+
m\log\frac{1+\frac{k}{2m}(1-a)+...}{1-\frac{k}{2m}(1+a)+...}\\
=&m\log\frac{1+a}{1-a}+k\log{2m}+
k\log\frac{2a}{1-a^2}+
k-k\log k+...
\end{split}
\end{equation}

Actually it was important for us to note that in the 
second (logarithmic) term in
asymptotic we have the factor $k$.

\bigskip

To finish this section let us discuss asymptotic for the case
$$
k=qm, \quad m\to \infty
$$
for a fixed $q$.
Note that now $c$ is just a constant
\begin{equation}
c=\sqrt{\frac{a^2+qa}{1+qa}}
\end{equation}
and we have
\begin{equation}\label{last}
B=m\left(q\log\frac{c+a}{c-a}+\log\frac{1+c}{1-c}\right),
\end{equation}
and $B=2m\log\frac{1+\sqrt{a}}{1-\sqrt{a}}$ for $q=1$.

%\section{Special addition}
%
%It looks like we are interested to find the  conformal
%map $\phi$ of the upper half plane onto the region in the upper
%half plane above the curve
%$$
%\gamma=\{u+iv=\phi(x):\ x\in\bbR\}
%$$
%such that
%$$
%L\sin v(x)\sinh u(x)=x,\quad x\in \bbR
%$$
%and normalized by $\phi(0)=0$, $\phi(z)\sim z$,
%$z\to\infty$.
%
%Let us rewrite the above equation in terms of the unknown function,
%say $\rho$ and
%its Hilbert transform $\tilde\rho$. We use the integral representation
%$$
%\phi(z)=z+\frac 1 \pi \int_0^\infty
%\left[\frac{1}{x-z}-\frac{1}{x+z}\right]v(x)\,dx.
%$$
%The curve  has the asymptote
%$v\to \pi$, $x\to\infty$. We define $\rho:=\pi-v$ to write
%$$
%\phi(z)=z+i\pi-\frac 1 \pi \int_0^\infty
%\frac{\rho\,dx}{x-z}.
%$$
%Finally since
%$$
%\frac 1 \pi \int_0^\infty
%\frac{\rho\,dx}{x-z}=-\tilde\rho+i\rho,
%$$
%we get $u(x)=\tilde\rho(x)+x$ and $v(x)=\pi-\rho(x)$. Thus the equation
%is
%$$
%L\sin\rho(x)\sinh(\tilde\rho(x)+x)=x.
%$$

\bibliographystyle{amsplain}

Address:\\[.1cm]
Abteilung f\"ur Dynamische Systeme \\
und Approximationstheorie\\
Universit\"at Linz\\
4040 Linz, Austria\\
{Franz.Peherstorfer@jku.at}

Address: \\[.1cm]
Abteilung f\"ur Dynamische Systeme \\
und Approximationstheorie\\
Universit\"at Linz\\
4040 Linz, Austria\\
and\\
Department of Mathematics \\
Bar Ilan University, Israel\\
Petro.Yudytskiy@jku.at

\end{document}